\documentclass{article}
\usepackage{amsmath, amssymb}
\usepackage{tikz-cd}
\usepackage{stmaryrd}
\usepackage{tikz}
\usetikzlibrary{shapes.geometric} 
\usepackage[a4paper,margin=2cm]{geometry}
\usepackage{hyperref}
\title{Some Structures arising from the Farey Fractal}
\author{Shai Haran}
\date{\today}

\begin{document}
 \maketitle

 \begin{abstract}

When one study the shape of Farey sequences, that is the shape of the binary planar rooted tree
\[
C(R) = \{ (x,y)\in \mathbb{N}\times \mathbb{N} \;\mid\; \gcd(x,y)=1,\; x+y\le R \},
\]
one rediscovers the Zeckendorf expansion of $R$ as a sum of non-consecutive Fibonacci numbers.\\
With $\varphi$ the golden ratio, we notice the subtle difference between identifying
\[
\varphi^2 \;\longleftrightarrow\; \varphi+1 , \quad \text{hence} \quad \varphi^{n+1} \;\longleftrightarrow\; \varphi^n + \varphi^{n-1},
\]
obtaining the totally disconnected, locally compact (pro finite) "Fibonadic numbers", \\
versus identifying also
\[
\varphi \;\longleftrightarrow\; \sum_{j\geq 0} \varphi^{-2j}, 
\quad \text{hence} \quad 
\varphi^{n+1} \;\longleftrightarrow\; \sum_{j\geq 0} \varphi^{\,n-2j},
\]
obtaining the continuum of non-negative real numbers.
 \\\\
\end{abstract}

Put
\[
\mathbb{N}^+ = \{1,2,3,\ldots\}, 
\qquad 
\mathbb{N}^{++} = \{2,3,4,\ldots\}.
\]

\medskip
\textbf{Idea:} instead of studying
\[
\mathbb{N}^{++} \;\backslash\; \mathbb{N}^{++} \cdot \mathbb{N}^{++} \;=\; \text{Primes},
\]

study
\[
\mathcal{X} = \mathbb{N}^+ \times \mathbb{N}^+ \;\backslash\; \mathbb{N}^{++} \cdot (\mathbb{N}^+ \times \mathbb{N}^+).
\]

That is,
\[
\mathcal{X} = \{\, (x,y)\in \mathbb{N}^+ \times \mathbb{N}^+ \;\mid\; \gcd(x,y)=1 \,\} 
\equiv \; \mathbb{Q}^+ 
\]
\medskip 
\[
(x,y) \;\longleftrightarrow\; \frac{y}{x}
\]

\medskip
\textbf{Thus $\mathcal{X}$ is the free abelian group on the set of primes, but it has more structure!}

\newpage
\section*{Fibonacci numbers}

\[
a_0 = 0, \; a_1 = 1, \; a_2 = 1, \; a_3 = 2, \; a_4 = 3, \; a_5 = 5
\]
\[
a_{n+1} = a_n + a_{n-1} ,\; all \; n
\]

\subsection*{Binet formula}
\[
a_n = \frac{\varphi^n - \psi^n}{\sqrt{5}}
\]
with $\varphi$ the golden ratio:
\[
\varphi = \frac{1+\sqrt{5}}{2}, \quad \varphi^2 = \varphi + 1
\]
and with $\psi$ its conjugate
\[
\quad \psi = \frac{1-\sqrt{5}}{2} \quad = \quad 1-\varphi = -\frac{1}{\varphi}
\]
 
Put $f_{-n}$ = 0 for all $n\geq 0$
\[
f_0 = 0, \; f_1 = 1, \; f_2 = 2, \; f_3 = 3, \; f_4 = 5, \; f_5 = 8, \; f_6 = 13
\]
\[
f_{n+1} = f_n + f_{n-1}, \quad n > 1
\]
\[
f_n = a_{n+1}, \quad n>0
\]

\subsection*{Zeckendorf expansion}
Let $R \in \mathbb{N}^+$.  
\[
M(R) = \max \{n \mid f_n \leq R\}
\]
\[
\delta R = R - f_{M(R)}
\]
\[
\Rightarrow \quad R = \sum_{j \ge 0} f_{M(\delta^{\,j} R)}
\]

Define
\[
Z(R) \in \{0,1\}^{\mathbb{N}^+}, \quad Z(R)_n = 1 \iff n = M(\delta^{\,j} R) \text{ for some } j .
\]
\[
R = \sum_{n \ge 1} Z(R)_n \, f_n = \sum_{n \ge 1} z_n f_n .
\]
Note: The sequence $ z_n$  has no consecutive 1's, $z_n \bullet z_{n+1} \equiv 0$ for all n
\subsection*{Definitions}
\[
\pi(R) = \sum_{n \ge 1} z_n f_{n-1}, 
\qquad 
j(R) = \sum_{n \ge 1} z_n f_{n+1}
\]
\[
\pi \circ j = id
\]
\[
Z(\pi(R))_n = Z(R)_{n+1}
\]
\[
Z(j(R))_n =
\begin{cases}
0 & n = 0,\\
Z(R)_{n-1} & n > 0
\end{cases}
\]

\bigskip
\noindent
\textbf{Thus $\pi$ erases the last digit of the Zeckendorf expansion of a number, and $j$ adds a zero to the expansion.}

\medskip
\noindent
\textbf{Note that the maps $\pi$ and $j$ do not preserve addition, but their coboundaries take values in $\{-1,0,1\}$.}

\subsection*{Lemma}
\[
\partial \pi(x,y) = \pi(x+y) - \pi(x) - \pi(y) \in \{-1,0,1\}
\]
\[
\partial j(x,y) = j(x+y) - j(x) - j(y) \in \{-1,0,1\}
\]

\newpage
\section*{Stern–Brocot Tree and Farey Graph}

\medskip
 
Let 
\[
SL_2(\mathbb{N}) = \left\{ 
\begin{pmatrix} 
a & b \\ 
c & d 
\end{pmatrix}
\,\middle|\, ad - bc = 1, \; a,b,c,d \in \mathbb{N} 
\right\}.
\]

It is the free monoid on two generators 
\[
\begin{pmatrix} 1 & 1 \\ 0 & 1 \end{pmatrix},
\quad 
\begin{pmatrix} 1 & 0 \\ 1 & 1 \end{pmatrix}.
\]

We have a bijection which we use as identification.

\[
SL_2(\mathbb{N})  
\quad \longleftrightarrow \quad 
X = \left\{ (x,y) \in \mathbb{N}^2 \;|\; \gcd(x,y)=1 \right\}
\]

\[
g \longmapsto (1,1)g,
\quad
 g_v = 
\begin{pmatrix} v^- \\ v^+ \end{pmatrix}
\;\longmapsfrom\; v
\]

\[
v^+ = (0,1) g_v
\in X \cup \{\infty = (0,1)\} \quad v^+ \quad \text{is the upper bound of $v$} 
\]

\[
v^- = (1,0) g_v
\in X \cup \{0 = (1,0)\} \quad v^- \quad \text{is the lower bound of $v$} 
\]

\[
v = v^+ + v^-, \quad \text{$v$ is the mediant of $v^+$ and $v^-$} 
\]

The planar binary tree structure, with vertices $X$ and root $(1,1)$,  
is obtained by viewing $v + v^+ > v + v^-$ as the descendants of $v$,  
and $v$ is the \textbf{\textit{Mother}} of $v + v^+$ and $v + v^-$.

\subsection*{Farey Graph}
\noindent
\textbf{The oriented graph's structure is given by viewing each $g \in X$ as an edge going from the lower bound $(1,0)g$ to the upper bound $(0,1)g$. 
}
\medskip
\medskip
\noindent
\textbf{It has edges $X$, and vertices}
\[
\overline{X} \;=\; X \,\sqcup\, \{\, 0 := (1,0), \; \infty := (0,1) \,\}.
\]
\medskip
\medskip
\noindent
\textbf{
Note that for $v \neq (1,1)$, among $v^+, v^-$ there is a longer vector, the \emph{Mother} of $v$ denoted $M(v)$, and a smaller one, called the \emph{Father} of $v$ and denoted $F(v)$.}

\subsection*{Finite subtrees}
Let $\mathcal{C} =$ finite subtrees of $X \quad \Longleftrightarrow \quad$ paths in $\bar{X}$ from $0$ to $\infty$.

\[
 c
\;\;\Longleftrightarrow\;\; 
\partial c = \{\, 0 = (1,0) \,\} \;\sqcup\; c \;\sqcup\; \{\, \infty = (0,1) \,\}.
\]
$c \in \mathcal{C}[n]$ if $\#c = n-1$,  $\#\partial c = n+1$, and $\partial c$ has $n$ steps.

\[
\mathcal{C} = \bigcup_{n \ge 1} \mathcal{C}[n]
\]
For $c \in \mathcal{C}$, let $\Phi c$ denote the set of branching points of $c$:
\[
\Phi c \;=\; \{\, v \in c \;\mid\; v+v^{+} \in c \ \text{and}\ v+v^{-} \in c \,\}.
\]
We concentrate on
\[
C(R) = \{ (x,y) \in X \;\mid\; |(x,y)| := x+y \leq R \}
\]
and on
\[
 C(R,n) \; := \; \Phi^n C(R) 
= \Bigl\{ v \in X \;\Big|\; |v|_n := a_{n+2}|M(v)| + a_{n+1}|F(v)| \leq R \Bigr\}.
\]

\medskip
\noindent
Note that the paths in the tree X outgoing from the root, that go fastest out are $\{(f_n, f_{n+1})\}$ and $\{(f_{n+1}, f_n)\}$ 
\newpage
\noindent
Thus:
\[
C(R, M(R)-1) = \varnothing, 
\qquad 
C(R, M(R)-2) = \{ (1,1) \}.
\]

\medskip
\noindent
There is an operad structure on $\mathcal{C}$, and
\[
C(R,n) = C(R,n+1)\cdot \lambda^{R}_{n}
\]
with “d.n.a”
\[
\lambda^{R}_{n} = \sum_{j} \lambda^{R}_{n,j}
\; \bullet \;
\bigl[c_j,c_{j+1}\bigr] \in \mathbb{Z} \; \partial  C(R,n+1).
\]
\noindent
That is, for $\bigl[c_j,c_{j+1}\bigr]\in \partial C(R,n+1)$ and $\lambda_j  =  \lambda^{R}_{n,j} \in\mathbb{Z}$ \\
\\
to construct $ \partial  C(R,n)$ from $ \partial  C(R,n+1)$ one replaces the step
\[
[c_j,c_{j+1}]
\]
by
\[
\begin{cases}
c_j \;\longrightarrow\; [c_j+c_{j+1}] \;\longrightarrow \cdots \longrightarrow [c_j+(\lambda_j+1)c_{j+1}] \longrightarrow c_{j+1}, & \lambda_j > 0, \\[6pt]
c_j \;\longrightarrow [c_j+c_{j+1}] \longrightarrow c_{j+1}, & \lambda_j=0 \quad \text{(just add mediant)}, \\[6pt]
c_j \;\longrightarrow [c_{j+1}+(|\lambda_j|+1)c_j] \;\longrightarrow \cdots \longrightarrow [c_{j+1}+c_j] \longrightarrow c_{j+1}, & \lambda_j < 0.
\end{cases}
\]

\subsection*{The $n$-th layer}
Define the $n$-th layer:
\[
\ell(n,R) := C(R, M(R)-n-1) \setminus C(R, M(R)-n)
\]
\[
= \{ v \in X \;\mid\; |v|_{M(R)-n-1} \leq R < |v|_{M(R)-n} \}.
\]
and the base
\[
b(n,R) = C(R, M(R) - n - 1)
       = \bigsqcup_{1 \leq m \leq n} \ell(m,R)
\]
\[
= \left\{\, v \in X \;\middle|\; 
   f_{M(R)-n} \; \big|M(v)\big| + f_{M(R)-n-1} \; |F(v)| \leq R
   \,\right\}.
\]

\medskip
\noindent
Note that for fixed $v$ the norms $|v|_k$ satisfy the same recurrence
\[
|v|_{k+1} = |v|_k + |v|_{k-1},
\]
hence
\[
\pi(|v|_{k+1}) = |v|_k,
\]
for $k$ large $\bigl(f_k > |v|^2\bigr)$, and since $\pi$ is order-preserving,
\[
R \in [\,|v|_k,\,|v|_{k+1}) \quad \text{if and only if} \quad \pi(R) \in [\,|v|_{k-1},\,|v|_k).
\]
Taking $k = M(R)-n$, we get for $M(R)\gg n$,
\[
\ell(n,R) = \ell(n,\pi(R)).
\]
\newpage
\section*{Fibonadic Numbers}

Define the (non–negative) \emph{Fibonadic Numbers} $\mathcal{F}$ as the inverse limit of $\mathbb{N}$ with respect to $\pi$:
\[
\mathcal{F} := \lim_{\longleftarrow} \bigl\{ \mathbb{N} \xleftarrow{\;\pi\;} \mathbb{N} \xleftarrow{\;\pi\;} \mathbb{N} \xleftarrow{\;\pi\;} \cdots \bigr\}.
\]

\medskip
\noindent
Alternatively, using Zeckendorf expansion:
\[
\mathcal{F} = \Bigl\{ z=(z_n) \in \{0,1\}^{\mathbb{Z}} \;\Big|\; 
z_n=0 \text{ for } n\gg 0,\; z_n z_{n-1}=0 \text{ for all } n \Bigr\}.
\]
For $z\ne z'$, let 
\[
\delta(z,z') := \max\{\,n \;\mid\; z_n \ne z'_n \,\}.
\] 
\\
\noindent
Then $d(z,z') := \varphi^{\delta(z,z')}$ is a metric, and $\mathcal{F}$ is complete.  \\
\\
Put $z < z'$ if for $n = \delta(z,z')$, $z_n=0$, $z'_n=1$. It's a total order on $\mathcal{F}$.\\

\noindent
For $z=(z_n)\in \mathcal{F}$, let $\{R^z_n\}_{n\in\mathbb{N}}$ be the associated $\pi$–coherent sequence:
\[
R^z_n := \sum_{j} z_j\, f_{n+j} \;\;\in \mathbb{N},
\]
and
\[
\pi(R^z_{n+1}) = R^z_n.
\]

\medskip
\noindent
We get a decomposition of $\mathcal{X}$ into levels:
\[
\mathcal{X} = \bigsqcup_{n\geq 1} \ell(n,z),
\]
with
\[
\ell(n,z) \equiv \ell(n,R^z_m)\quad \text{for any $m\gg n$}.
\]
and base 
\[
b(n,z) = \bigsqcup_{1\leq m \leq n} \ell(m,z)
\]

\medskip
We let the levels function
\[
f^z : \mathcal{X} \;\longrightarrow\; \mathbb{N}
\]
be defined by
\[
f^z(v) = n \quad \Longleftrightarrow \quad v \in \ell(n,z).
\]
Thus
\[
b(n,z)=\{v\; \Big| \; |v|_{M(R^z_m)-n} \leq R^z_m \; for \; m \gg n\}
\]
Note that
\[
f^z(v+M(v)) = f^z(v) + 1
\]
and
\[
f^z(v+F(v)) = f^z(v) + d f^z(v), 
\qquad d f^z(v) \in \{0,1\}.
\]
\\
\\
Thus $f^z$ is determined by $d f^z$ taking values $0,1$, or by the subset
\[
Z(z) = \{\, v \in \mathcal{X} \;\mid\; d f^z(v) = 0 \,\}.
\]
Moreover,
\[
Z(z) \cap \ell(n,z) 
= \bigl\{\, v \;\big|\; |v+F(v)|_{M(R^z_m)-n} \;\leq\; R^z_m \;<\; |v|_{M(R^z_m)-n+1},\; m \gg n \,\bigr\}.
\]

\newpage
\subsection*{Normalization Lemma}

\medskip
Let $\mathbb{N}((t))$ denote the commutative rig (i.e. ring without negatives) of Laurent series with coefficients in $\mathbb{N}$. \\
 \\
Let

\[
\mathcal{B} := \Bigl\{ z(t)=\sum_n z_n t^n \in \mathbb{N}((t))\ \Big|\ 
V_\varphi(z) := \sum_n z_n \varphi^n < \infty \Bigr\}.
\]
This is a sub–rig of $\mathbb{N}((t))$, and we have a surjective rig homomorphism
\[
V_\varphi : \mathcal{B} \;\longrightarrow\; \mathbb{R}_{\ge 0} := [0,\infty).
\]
We can view $z=(z_n)\in\mathcal{F}$ as a power series $z(t)=\sum_n z_n t^n,$\\
\\
making $\mathcal{F}$ a subset of $\mathcal{B}$, and the $\varphi$–value is finite: $for \; z \in F, \;V_\varphi(z)=\sum_n z_n \varphi^n < \infty.
$\\
\\
We have a canonical projection 
\[
P: \mathcal{B} \longrightarrow \mathcal{F},
\]
preserving the $\varphi$–value:
\[
V_\varphi(z)=V_\varphi(P(z)).
\]
For $z = (z_n) \in \mathbb{N}^{\mathbb{Z}}$, with 
\[
\sum_n z_n \varphi^n < \infty,
\]
(hence $z_n = 0$ for $n \gg 0$), making the moves

\begin{itemize}
  \item[(A)] 
  \[
  (\dots 0\, z_n \, z_{n-1}\, z_{n-2} \dots) 
  \;\;\longrightarrow\;\;
  (\dots 1\,  (z_n-1)\, (z_{n-1}-1)\, z_{n-2}\dots)
  \]
  if $z_n \geq 1$, $z_{n-1} \geq 1$,
  
  \item[(B)]
  \[
  (\dots 0\, z_n \,0\, z_{n-2} \dots) 
  \;\;\longrightarrow\;\;
  (\dots 1\, (z_n-2)\, 0\, (z_{n-2}+1)\dots)
  \]
  if $z_n \geq 2$,
\end{itemize}
produces a sequence that converges to $\mathcal{P}(z) \in \mathcal{F}$.

\medskip
\noindent
Note that steps (A), (B) preserve the $\varphi$-value. 
\\\\\\
Say that the sequence $z$ is OK above $z_n$ if 
\[z_m \in \{0,1\}, \; \text{and}\; z_m \;\bullet \; z_{m-1}=0 \; \text{ for all} \;m > n  
\]
If we are in the case $z_n \geq 1$ and $z_{n-1} \geq 1$, a finite number of $A$ moves will produce a sequence that is OK above $z_{n-1}$, or a sequence OK above $z_n$, with $z_n \geq 2$, $z_{n-1} = 0$. The last case is when move $B$ applies: 
it will replace $z_n$ by $z_n-2$, and after a few cleanup $A$-moves it will be OK above $z_n$. Repeating this we see that given $(z)$ OK above $z_n$, 
it will be transformed to a sequence OK above $z_{n-1}$.
\\\\
Note that for fixed $m \in \mathbb{Z}$, as $n \to -\infty$ the values of our sequence above $m$ 
can change, in which case the sum 
\[
\sum_{j \geq m} z_j \varphi^j
\]
only increases, by an amount $\geq \varphi^m$. 
\newpage
\noindent
The boundedness of the $\varphi$-value 
\[
\sum_n z_n \varphi^n < \infty
\]
guarantees that this canonical process stabilizes above any finite index $m$. 
\\\\
Denoting by $P_n(z)$ the sequence obtained from $z$ after making it OK above $z_n$ in a fixed stable way, i.e. the $z_m$ for $m > n$ get their final value, so 
\[
P_n(z) \xrightarrow[n \to -\infty]{} P(z)
\]
and for any finite index $m$, the coefficient
\[
P(z)_m = P_n(z)_m \quad \text{for all } n <m
\]

\subsection*{Rig structure on $\mathcal{F}$}

Define for $z, z' \in \mathcal{F}$:
\[
z+z' := P\!\bigl(z(t)+z'(t)\bigr),
\qquad 
z\cdot z' := P\!\bigl(z(t)\cdot z'(t)\bigr).
\]

These operations make $\mathcal{F}$ into a commutative rig, and $P$ is a rig homomorphism factoring $V_\varphi$:

\[
\begin{tikzcd}[row sep=3em, column sep=4em]
\mathcal{B} \arrow[r,"V_\varphi"] \arrow[d,"P"'] & \mathbb{R}_{\ge 0} \\
\mathcal{F} \arrow[ur,"V_\varphi"'] &
\end{tikzcd}
\]

\medskip
The operations $+, \bullet$ are clearly commutative, with unit for addition 
\[
0 = (z_n \equiv 0) \quad \text{the zero sequence,}
\] 
and unit for multiplication 
\[
1 = (z_n \equiv \delta_n) \quad \text{the sequence with } z_0 = 1, \; z_n = 0 \; \text{for } n \neq 0.
\]
To see the associativity of these operations, and the distributive law, one uses the identity
\[
P(x \circ y) = P(x \circ P(y)), \quad x,y \in B,
\]
with $\circ$ either addition or multiplication. 
\\
Notice that if $z'$ is obtained from $z$ by a finite number of moves, each exchanging 
\[
t^{n+1} \leftrightarrow t^n + t^{n-1},
\]
we have $P(z') = P(z)$, and $P_j(z') = P_j(z)$ for $j \ll 0$.
\\
Thus
\[
P(x \circ y) = P(x \circ P_n(y)), \quad \text{all } n,
\]
but for fixed $m$, the final value of the $m$-coefficient involve only finitely many of the other coefficients, and so 
\[
P(x \circ P_n(y))_m = P(x \circ P(y))_m
\]
for $n \ll N(m, V_\varphi(y), V_\varphi(x)).$
\newpage
\noindent
Let $\mathbb{N}[t,t^{-1}] \subseteq \mathcal{B}$ be the sub–rig generated by $t$ and $t^{-1}$. \par
\vspace{1\baselineskip}
\medskip
\noindent
Then $\mathcal{N} := P(\mathbb{N}[t,t^{-1}])$ is the sub–rig of $\mathcal{F}$ consisting of finite sequences:
\[
\mathcal{N} = \Bigl\{ z=(z_n)\in \mathcal{F} \;\Big|\; z_n=0,\ z_{-n}=0 \ \text{for } n\gg 0 \Bigr\}.
\]

\medskip
\noindent
For $z\in \mathcal{N}\setminus\{0\}$, let $n_0=\delta(z)=\min\{n \;\mid\; z_n=1\}$, and define $z^-=(z_n^-)\in \mathcal{F}$ by
\[
z^-_n =
\begin{cases}
z_n, & n>n_0, \\[6pt]
0, & n\leq n_0,\; n\equiv n_0 \pmod{2}, \\[6pt]
1, & n\leq n_0,\; n\not\equiv n_0 \pmod{2}.
\end{cases}
\]

\medskip
\noindent
Note that $z^- < z$, but they have the same $\varphi$–value:
\[
V_\varphi(z) = V_\varphi(z^-).
\]
The equivalence relation on $\mathcal{F}$ generated by
\[
z \sim z^- \qquad \text{for } z \in \mathcal{N}
\]
is compatible with the operations of addition and multiplication, and we have an identification
\[
\mathcal{F} \;\;\longrightarrow\;\; \mathcal{F}/\!\sim \;\; = \;\;\mathbb{R}_{\ge 0}
\]

\textbf{Beware:} Operations are \emph{not} continuous on $\mathcal{F}$.

\medskip
\noindent
E.g.
\[
\lim_{n \to \infty} \bigl(1^- + \Phi^{-n}\bigr) 
= 1 \;>\; 1^- = 1^- + 0 
= 1^- + \lim_{n \to \infty} \Phi^{-n}.
\]
\\
\subsection*{Quotients and Fundamental Domain}
\medskip
Let $\Phi \in \mathcal{F}$ be the image of $t$, so that $\Phi^m \cdot \Phi^n = \Phi^{m+n}$, and for $z \in \mathcal{F}$, $\Phi^n z$ is obtained by shifting its digits by $n$. \\\\
Note that 
\[
R^{\Phi z}_m = R^z_{m+1} \;,
b(n,\Phi z) = b(n,z)\;, f^{\Phi z} = f^z.
\]
Define:
\[
S := \bigl(\mathcal{F} \setminus \{0\}\bigr) \big/ \{\Phi^n\}_{n \in \mathbb{Z}}.
\]
Alternatively, let the "Principal Units"
\[
\mathcal{D} := \{\, z=(z_n)\in \mathcal{F} \;\mid\; z_0=1 \;\;\text{and}\;\; z_n=0 \;\;\text{for } n>0 \,\}.
\]
It is a fundamental domain for the action of $\{\Phi^n\}_{n\in\mathbb{Z}}$ on $\mathcal{F}\setminus\{0\}$, and we have an identification
\[
\mathcal{D} \;\;\longleftrightarrow\;\; S.
\]
We have multiplication on $S$, and taking $\varphi$–values gives a homomorphism onto the circle:
\[
V_\varphi : S \;\longrightarrow\; \mathbb{R}^+ \big/ \{\varphi^n\}_{n\in\mathbb{Z}}.
\]

\medskip
\noindent
We have a total order on $\mathcal{D}$, with minimal element $1=\Phi^0$, and maximal element
\[
(\Phi^1)^- = (\ldots 0\ldots 01010101\ldots 01\ldots),
\]
and taking $\varphi$–values gives an order preserving map
\[
V_\varphi : \mathcal{D} \;\longrightarrow\; [1,\varphi].
\]
For $z \in \mathcal{D}$, we have
\[
\qquad M(R_m^z) \equiv m 
\]
and the base
\[
b(n,1) \;\subseteq\; b(n,z) \;\subseteq\; b\!\left(n,\Phi^-\right) 
\]
Indeed for $z, z' \in \mathcal{D} $
\[
z \leq z' \; \;\text{implies} \;\; b(n,z) \subseteq b(n,z'),\;  and \; f^z(v) \geq f^{z'}(v) \; \text{for all} \; v 
\]
\\\\
\noindent
\[
b(n,z)
= \left\{ v \ \middle|\  
f_{m-n+1}\cdot |M(v)| + f_{m-n}\cdot |F(v)|
\le \sum_j f_{m+j}\cdot z_j \quad \text{for } m \ge n
\right\}
\]

\[
= \left\{ v \ \middle|\ 
\varphi^{m-n+2}\cdot |M(v)| + \varphi^{m-n+1}\cdot |F(v)|
\le \varphi^
{m+1} \sum_j \varphi^j z_j
\quad \text{for all } m
\right\}
\]

\[
= \left\{ v \ \middle| \ 
|v|_\varphi := |M(v)| + \varphi^{-1}|F(v)|
\le \varphi^{n-1}\cdot V_\varphi(z) \right\}
\]

with the $\varphi$-norm

\[
|v|_\varphi = |M(v)| + \varphi^{-1}|F(v)|
\equiv |M(v)-F(v)| + |F(v)|\cdot \varphi
\in N[\varphi] \subsetneq N[\varphi,\varphi^{-1}]
\]

\[
\ell(n,z)
= \left\{ v \ \middle|\ 
|v|_\varphi \le \varphi^{n-1}\cdot V_\varphi(z)
< \varphi\cdot |v|_\varphi \right\}
\]

\[
Z(z)\cap \ell(n,z)
= \left\{ v \ \middle|\ 
|v|_\varphi + |F(v)|
\le \varphi^{n-1}\cdot V_\varphi(z)
< \varphi\cdot |v|_\varphi \right\}
\]
This is the precise way in which the base-$\varphi$-expansion of Bergman~[2],
(which we used for defining the $\varphi$-value $V_{\varphi}(z)$, and for the
arithmetic of the rig of Fibonadic numbers), gives the asymptotics of the
Zeckendorf expansion. \newline

\noindent
The element \( v \in X \) is relevant for the \(n\)-layer if and only if

\[
|v|_\varphi \in [\varphi^{n-1},\, \varphi^n].
\]

\vspace{1\baselineskip}
\medskip
\noindent
Let $B(n)$ denote the set of possible configurations of height $n$,
\[
B(n) = \{\, b(n,z) \mid z \in \mathcal{D} \,\}.
\]
We have restriction map $r : B(n) \to B(n-1)$, making 
\[
\mathbb{B} = \bigcup_{n} B(n)
\]
into a planar rooted tree,

\[
\text{root: } B(1) = \{\, b(1,z) = \{(1,1)\} \;\; \text{any } z \in \mathcal{D} \,\}.
\]
Order on $B(n)$: 
\[
b \leq b' \;\;\Longleftrightarrow\;\; r^j b \;\subseteq\; r^j b', \quad j=0,\ldots,n.
\]
\\\\
The set $\mathcal{D}$ is the set of paths outgoing from the root of a 
planar binary rooted tree with vertices the finite $0,1$ expansions, 
with root $10$, and for a finite expansion $z$ we have its descendants $z0$ and $z10$.

The map 
\[
\mathcal{D} \longrightarrow B(n), 
\quad z \longmapsto b(n,z)
\]
is continuous and order preserving. 
\\\\
It follows that there are unique $z^n_j \in \mathcal{N} = \mathbb{N}[\varphi, \varphi^{-1}]$
\[
1 = z^n_0 < z^n_1 < \cdots < z^n_k < z^n_{k+1} < \cdots < z^n_{\#B(n)} = \Phi
\]
such that for $z \in \mathcal{D}$:
\[
z^n_k \leq  z \leq (z^n_{k+1})^- 
\;\;\Longleftrightarrow\;\; 
b(n,z^n_k) = b(n,z) = b\!\left(n,(z^n_{k+1})^-\right).
\]
\\
Note that $\{ z^{n-1}_j \}$ is a subset of $\{ z^n_k \}$, 
and for $z^n_k \leq z \leq (z^n_{k+1})^-$
\[
r\,b(n,z) = b(n-1,z) = b\!\left(n-1,z^{n-1}_j\right)
\quad \iff \quad
[z^n_k , (z^n_{k+1})^-] \subseteq [z^{n-1}_j, (z^{n-1}_{j+1})^-].
\]
\\
The new \( z_k^n \)'s are given as the \(\varphi\)-norms
\[
|v|_{\varphi} \in \mathbb{N}[\varphi]
\;\subsetneq\;
\mathbb{N}[\varphi,\varphi^{-1}]
\;\subsetneq\;
\mathbb{Z}[\varphi]
\;\subseteq\;
\mathbb{R}\times\mathbb{R},
\qquad
p(\varphi) \longmapsto \bigl(p(\varphi),\, p(\psi)\bigr),
\]
for \( v \in X \) which is relevant for the \(n\)-layer.
\vspace{1\baselineskip}
\newline
\noindent
Let $\varprojlim \mathbb{B}$ denote the set of paths in $\mathbb{B}$ outgoing from the root,  
it is the set of configurations of infinite height, or equivalently  
the set of functions 
\[
f : X \to \mathbb{N},
\]
\[
f(v) = \min \{\, n \mid v \in b(n) \,\}
\]
\\
The map
\[
\begin{aligned}
\mathcal{D} &\;\longrightarrow\; \varprojlim \mathbb{B}, \\
z           &\;\longmapsto\; \{ b(n,z) \}_n \;\longleftrightarrow\; f^z
\end{aligned}
\]
is a bijection.
\newpage
\subsection*{The Tree of configurations}
We give the beginning of the tree $\mathbb{B}$ with the values of $z^n_j$'s, for $n \;=\; 1,2,3,$ and with a \tikz[baseline={(0,-0.5ex)}]{\node[circle,draw,double,minimum size=6pt,inner sep=0pt] (c) {};} denoting points of $Z(z)$ that determine the configuration, that is the points $v$ such that $v \; +\; F(v)$ is in the same level as $v$, exhibiting the linear growth of $Z(z)$.

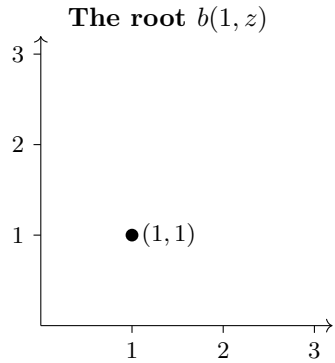
\begin{center}
\textbf{The root $b(1,z)$}

\begin{tikzpicture}[scale=1.2]
  \draw[->] (0,0)--(3.2,0); \draw[->] (0,0)--(0,3.2);
  \foreach \x in {1,2,3} \draw (\x,0) -- (\x,-.08) node[below] {\small \x};
  \foreach \y in {1,2,3} \draw (0,\y) -- (-.08,\y) node[left] {\small \y};
  \fill (1,1) circle (2pt) node[right] {\small $(1,1)$};

\end{tikzpicture}
 
\textit{It branches to two in the first level}
\end{center}
\vspace{1ex}
\begin{center}
$b(2,z)$
\end{center}
\begin{minipage}{0.48\linewidth}
\centering
\[
1 = z^{2}_{0} \leq z < (z^{2}_{1})^{-} = 10010\ldots 10\ldots
\]
\[
\text{i.e.} \;z=100\ldots 
\]
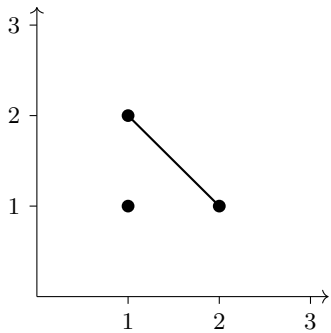
\begin{tikzpicture}[scale=1.2]
  \draw[->] (0,0)--(3.2,0); \draw[->] (0,0)--(0,3.2);
  \foreach \x in {1,2,3} \draw (\x,0) -- (\x,-.08) node[below] {\small \x};
  \foreach \y in {1,2,3} \draw (0,\y) -- (-.08,\y) node[left] {\small \y};
  \fill (1,1) circle (2pt);
  \fill (1,2) circle (2pt);
  \fill (2,1) circle (2pt);
  \draw[thick] (1,2) -- (2,1);
\end{tikzpicture}
\end{minipage}
\hfill
\begin{minipage}{0.48\linewidth}
\centering
\[
10100\ldots 0\ldots = z^{2}_{1} \leq z \leq z^{2}_{2} = \Phi^{-}
\]
\[
\text{i.e.} \;z=1010\ldots
\]
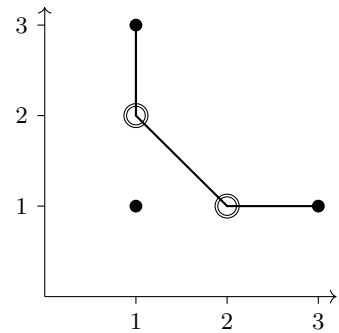
\begin{tikzpicture}[scale=1.2]
  \draw[->] (0,0)--(3.2,0);
  \draw[->] (0,0)--(0,3.2);
  \foreach \x in {1,2,3} \draw (\x,0) -- (\x,-.08) node[below] {\small \x};
  \foreach \y in {1,2,3} \draw (0,\y)--(-.08,\y) node[left] {\small \y};

  \fill (1,3) circle (2pt);
  \fill (3,1) circle (2pt);
   \fill (1,1) circle (2pt);

  \node[circle,draw,double,minimum size=8pt,inner sep=0pt] at (1,2) {};
  \node[circle,draw,double,minimum size=8pt,inner sep=0pt] at (2,1) {};

  \draw[thick] (1,3)--(1,2)--(2,1)--(3,1);
\end{tikzpicture}
\end{minipage}

\vspace{1ex}
\begin{center}
\textit{Each of these, branches to three in the next level.}
\end{center}

\newpage
\newgeometry{top=1.5cm, bottom=1.5cm}
\begin{center}
\textbf {$b(3,z)$}
\end{center}
\begin{minipage}{0.48\linewidth}
\centering
\[
 z=100\ldots
\]
\begin{tikzpicture}[scale=1.2]
  \draw[->] (0,0)--(3.6,0);
  \draw[->] (0,0)--(0,3.6);

  \foreach \x in {1,2,3} \draw (\x,0)--(\x,-.08) node[below] {\small \x};
  \foreach \y in {1,2,3} \draw (0,\y)--(-.08,\y) node[left]  {\small \y};

  \fill (1,3) circle (2pt);
  \fill (2,3) circle (2pt);
  \fill (3,2) circle (2pt);
  \fill (3,1) circle (2pt);
  \draw[very thick] (1,3)--(2,3)--(3,2)--(3,1);

  \fill (1,2) circle (2pt);
  \fill (2,1) circle (2pt);
  \draw[very thick] (1,2)--(2,1);

  \fill (1,1) circle (2pt);
\end{tikzpicture}
\vspace{0.3em}
\[
1=z^{3}_{0} \leq z
\]
\[
z \leq (z^{3}_{1})^-= 
10000010
\ldots 10\ldots
\]
\vspace{1em}
\begin{center}
\begin{tikzpicture}[scale=1]
  \draw[->] (0,0)--(4.6,0);
  \draw[->] (0,0)--(0,4.6);

  \foreach \x in {1,2,3,4} \draw (\x,0)--(\x,-.08) node[below] {\small \x};
  \foreach \y in {1,2,3,4} \draw (0,\y)--(-.08,\y) node[left]  {\small \y};

  \draw[very thick] (2,1)--(1,2);
  \foreach \p in {(2,1),(1,2)} \fill \p circle (2pt);

  \draw[very thick] (4,1)--(3,1)--(3,2)--(2,3)--(1,3)--(1,4);

  \foreach \q in {(4,1),(3,2),(2,3),(1,4)} \fill \q circle (2pt);
 \fill (1,1) circle (2pt);
  \node[circle,draw,double,minimum size=8pt,inner sep=0pt] at (3,1) {};
  \node[circle,draw,double,minimum size=8pt,inner sep=0pt] at (1,3) {};
\end{tikzpicture}\vspace{0.3em}
\[
1000010\ldots 0\ldots \;=\; z^{3}_{1} 
\;\leq\; z \]
\[z \leq\; (z^{3}_{2})^{-} 
\;=\; 10010010\ldots10\ldots
\]
\end{center}
\vspace{1em}
\begin{center}
\begin{tikzpicture}[scale=0.8]
  \draw[->] (0,0)--(5.6,0);
  \draw[->] (0,0)--(0,5.6);

  \foreach \x in {1,2,3,4,5} \draw (\x,0)--(\x,-.08) node[below] {\small \x};
  \foreach \y in {1,2,3,4,5} \draw (0,\y)--(-.08,\y) node[left]  {\small \y};

  \draw[very thick] (2,1)--(1,2);
  \foreach \p in {(2,1),(1,2)} \fill \p circle (2pt);

  \draw[very thick] (5,1)--(4,1)--(3,1)--(3,2)--(2,3)--(1,3)--(1,4)--(1,5);

  \foreach \q in {(5,1),(3,2),(2,3),(1,5)} \fill \q circle (2pt);
 \fill (1,1) circle (2pt);
  \node[circle,draw,double,minimum size=8pt,inner sep=0pt] at (4,1) {};
  \node[circle,draw,double,minimum size=8pt,inner sep=0pt] at (3,1) {};
  \node[circle,draw,double,minimum size=8pt,inner sep=0pt] at (1,3) {};
  \node[circle,draw,double,minimum size=8pt,inner sep=0pt] at (1,4) {};
\end{tikzpicture}\vspace{0.2em}
\[
10010100\ldots 0 \ldots \;=\; z^{3}_{2} \;\leq\; z \]
\[z \; \leq\; (z^{3}_{3})^{-} =\; 10010\ldots 10\ldots
\]
\end{center}
\end{minipage}
\hfill
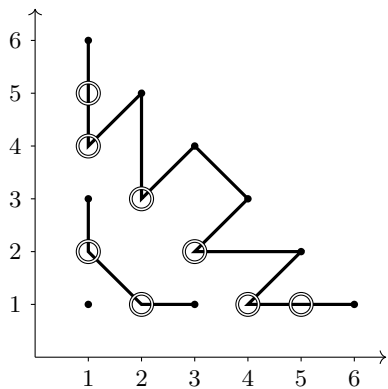
\begin{minipage}{0.48\linewidth}
\centering
\[
 z=101\ldots
\]
\begin{tikzpicture}[scale=0.8]
  \draw[->] (0,0)--(5.6,0);
  \draw[->] (0,0)--(0,5.6);

  \foreach \x in {1,2,3,4,5} \draw (\x,0)--(\x,-.08) node[below] {\small \x};
  \foreach \y in {1,2,3,4,5} \draw (0,\y)--(-.08,\y) node[left]  {\small \y};

  \draw[very thick] (3,1)--(2,1)--(1,2)--(1,3);

  \foreach \p in {(1,1),(3,1),(1,3)} \fill \p circle (2pt);

  \node[circle,draw,double,minimum size=8pt,inner sep=0pt] at (2,1) {};
  \node[circle,draw,double,minimum size=8pt,inner sep=0pt] at (1,2) {};

  \draw[very thick] (5,1)--(4,1)--(5,2)--(3,2)--(2,3)--(2,5)--(1,4)--(1,5);

  \foreach \q in {(5,1),(5,2),(3,2),(2,3),(2,5),(1,5)} \fill \q circle (2pt);

  \node[circle,draw,double,minimum size=8pt,inner sep=0pt] at (4,1) {};
  \node[circle,draw,double,minimum size=8pt,inner sep=0pt] at (1,4) {};
\end{tikzpicture}\vspace{0.3em}
\[
\text{1010}\ldots 0\ldots \;=\; z^{3}_{3} \;\leq\; z \]
\[z \leq\; (z^{3}_{4})^{-} \;=\; 10100010\ldots 10\ldots
\]
\vspace{1em}

\begin{center}
\begin{tikzpicture}[scale=0.8]
  \draw[->] (0,0)--(5.6,0);
  \draw[->] (0,0)--(0,5.6);

  \foreach \x in {1,2,3,4,5} \draw (\x,0)--(\x,-.08) node[below] {\small \x};
  \foreach \y in {1,2,3,4,5} \draw (0,\y)--(-.08,\y) node[left]  {\small \y};

  \draw[very thick] (3,1)--(2,1)--(1,2)--(1,3);

  \foreach \p in {(1,1),(3,1),(1,3)} \fill \p circle (2pt);

  \node[circle,draw,double,minimum size=8pt,inner sep=0pt] at (2,1) {};
  \node[circle,draw,double,minimum size=8pt,inner sep=0pt] at (1,2) {};

  \draw[very thick]
    (5,1)--(4,1)--(5,2)--(3,2)--(4,3)--(3,4)--(2,3)--(2,5)--(1,4)--(1,5);

  \foreach \q in {(5,1),(5,2),(4,3),(3,4),(2,5),(1,5)} \fill \q circle (2pt);

  \node[circle,draw,double,minimum size=8pt,inner sep=0pt] at (4,1) {};
  \node[circle,draw,double,minimum size=8pt,inner sep=0pt] at (3,2) {};
  \node[circle,draw,double,minimum size=8pt,inner sep=0pt] at (2,3) {};
  \node[circle,draw,double,minimum size=8pt,inner sep=0pt] at (1,4) {};
\end{tikzpicture}\vspace{0.2em}
\[
1010010\ldots 0\ldots \;=\; z^{3}_{4}
\;\leq\; z \]
\[z \leq\; (z^{3}_{5})^{-}
\;=\; 10101000\,10\ldots 10\ldots
\]
\end{center}
\vspace{1em}
\begin{center}
\begin{tikzpicture}[scale=0.7]
  \draw[->] (0,0)--(6.6,0);
  \draw[->] (0,0)--(0,6.6);

  \foreach \x in {1,2,3,4,5,6} \draw (\x,0)--(\x,-.08) node[below] {\small \x};
  \foreach \y in {1,2,3,4,5,6} \draw (0,\y)--(-.08,\y) node[left]  {\small \y};

  \draw[very thick] (3,1)--(2,1)--(1,2)--(1,3);

  \foreach \p in {(1,1),(3,1),(1,3)} \fill \p circle (2pt);

  \node[circle,draw,double,minimum size=8pt,inner sep=0pt] at (2,1) {};
  \node[circle,draw,double,minimum size=8pt,inner sep=0pt] at (1,2) {};

  \draw[very thick]
    (6,1)--(5,1)--(4,1)--(5,2)--(3,2)--(4,3)--(3,4)--(2,3)--(2,5)--(1,4)--(1,5)--(1,6);

  \foreach \q in {(6,1),(5,2),(4,3),(3,4),(2,5),(1,6)} \fill \q circle (2pt);

  \node[circle,draw,double,minimum size=8pt,inner sep=0pt] at (5,1) {};
  \node[circle,draw,double,minimum size=8pt,inner sep=0pt] at (4,1) {};
  \node[circle,draw,double,minimum size=8pt,inner sep=0pt] at (3,2) {};
  \node[circle,draw,double,minimum size=8pt,inner sep=0pt] at (2,3) {};
  \node[circle,draw,double,minimum size=8pt,inner sep=0pt] at (1,4) {};
  \node[circle,draw,double,minimum size=8pt,inner sep=0pt] at (1,5) {};
\end{tikzpicture}\vspace{0.2em}
\[
10101001000\ldots0\ldots \;=\; z^{3}_{5}
\;\leq\; z \]
\[z \leq\; (z^{3}_{6})^{-}
\;=\; \Phi^{-} \;=\; 1010\ldots 10\ldots
\]
\end{center}
 
\end{minipage}
\restoregeometry 

\end{document}